\documentclass[11pt,leqno]{article} 
\usepackage{graphics}
\newtheorem{thm}{Theorem}[section]
\newtheorem{lma}{Lemma}[section]

\newcommand{\beqa}{\begin{eqnarray}}
\newcommand{\eeqa}{\end{eqnarray}}

\newcommand{\pf}{\noindent {\bf Proof:} $\s$ }
\newcommand{\epf}{ \hfill$\diamondsuit$ \medskip}

\newcommand{\beq}{\begin{equation}}
\newcommand{\eeq}{\end{equation}}
\newcommand{\lbl}{\label}
\newcommand{\s}{\; \;}

\newcommand{\ep}{\epsilon}

\newcommand{\la}{\lambda}
\newcommand{\mb}{\mbox}
\newcommand{\ra}{\rightarrow}
\newcommand{\al}{\alpha}
\newcommand{\p}{\varphi}
\title{Infinitely many solutions for three classes of self-similar equations, with the $p$-Laplace operator}

\author{
Philip Korman   \\ 
Department of Mathematical Sciences \\ 
University of Cincinnati \\ 
Cincinnati Ohio 45221-0025 \\
kormanp@ucmail.uc.edu
}

\date{}

\begin{document}

\maketitle
\begin{abstract} 
We study the global solution curves, and prove the existence of infinitely many positive solutions for three classes of self-similar equations, with $p$-Laplace operator. In case $p=2$, these are  well-known problems involving the Gelfand equation,  the equation modeling electrostatic micro-electromechanical systems (MEMS), and a  polynomial nonlinearity.  We extend  the classical results of D.D. Joseph and T.S. Lundgren \cite{JL}  to the case $p \ne 2$, and  we generalize    the main result of Z. Guo and J. Wei \cite{GW} on the equation modeling MEMS.
 \end{abstract}

\begin{flushleft}
Key words:  Parameterization of the global solution curves, infinitely many solutions. 
\end{flushleft}

\begin{flushleft}
AMS subject classification: 35J60, 35B40.
\end{flushleft}

\section{Introduction}

We consider radial solutions on a ball in $R^n$ for three special classes of equations, involving the $p$-Laplace operator, the ones self-similar under scaling.
We now explain our approach for one of the  classes,  involving  the $p$-Laplace version of the equation which arises in modeling of electrostatic micro-electromechanical systems (MEMS), see \cite{P}, \cite{GG}, \cite{GW}  (with $p>1$, $\al>0$, $q>0$, $u=u(x)$, $x \in R^n$, $n \geq 1$)
\beq
\lbl{0.1}
\s\s\s \mbox{div} \left( |\nabla u|^{p-2} \nabla u \right)+\la \frac{|x|^{\alpha}}{(1-u)^q}=0, \; \mbox{for $|x|<1$} \s u=0, \; \mbox{when $|x|=1$} \,.
\eeq
Here $\la$ is a positive parameter. We are looking for solutions satisfying  $0<u<1$.
Radial solutions of this equation satisfy
\beq
\lbl{0.2}
\p (u')'+\frac{n-1}{r}\p (u') + \la \, \frac{r^{\alpha}}{(1-u)^q} =0 \s\s \mb{for $0<r<1$} \,, 
\eeq
\[
u'(0)=u(1)=0 \,, \s 0<u(r)<1  \,,
\]
with $\p (v)=v |v|^{p-2}$. It is easy to see that $u'(r)<0$ for all $0<r<1$,  which implies that the value of $u(0)$ gives the maximum value (or the $L^{\infty}$ norm) of our solution. Moreover, $u(0)$ is a {\em global parameter}, i.e., it uniquely identifies the solution pair $(\la ,u(r))$, see e.g., P. Korman \cite{K1}.
It follows that a two-dimensional curve in the $(\la,u(0))$ plane completely describes the solution set of (\ref{0.2}). The self-similarity of this equation allows one to parameterize the global solution curve, using the solution of a single initial value problem:
\beq
\lbl{0.3}
\p (w')'+\frac{n-1}{t}\p (w')=  \frac{t^{\al}}{w^q}, \s w(0)=1, \s w'(0)=0 \,.
\eeq
Its solution $w(t)$ is a positive and increasing function, which can be easily computed numerically. Following J.A. Pelesko \cite{P}, we show that the global solution curve of (\ref{0.2}) is given by
\[
(\la,u(0))=\left(\frac{t ^{\al +p}}{w^{p+q-1}(t)}\,, 1-\frac{1}{w(t)} \right) \,,
\]
parameterized by  $t \in (0,\infty)$.
In particular, $\la=\la(t)=\frac{t ^{\al +p}}{w^{p+q-1}(t)}$, and 
\[
\la '(t)=t^{\al +p-1}w^{-p-q} \left[(\al +p)w-t(p+q-1)w' \right] \,,
\]
so that the solution curve travels to the right (left) in the $(\la,u(0))$ plane if  $(\al +p)w-t(p+q-1)w'>0$ ($<0$). This makes us  interested in the roots of the function $(\al +p)w-t(p+q-1)w'$. If we set this function to zero
\[
(\al +p)w-t(p+q-1)w'=0 \,,
\]
then the general solution of this equation is 
\[
w(t)=c t^{\beta }, \s\s \beta=\frac{ \al +p}{p+q-1} \,.
\]
Quite remarkably, if we choose the constant $c=c_0=\left[\frac{1}{ \beta ^{p-1} \left[(p-1)(\beta-1)+n-1 \right]} \right]^{\frac{1}{p+q-1}}$
then 
\[
w_0(t)=c_0 t^{\beta }
\]
also solves the equation in (\ref{0.3}), along with $w(t)$. We show that $w(t)$ tends to $w_0(t)$ as $t \ra \infty$, and the solution curve of (\ref{0.2}) makes infinitely many turns if and only if $w(t)$ and $w_0(t)$ intersect infinitely many times. We give a sharp condition for that to happen, thus generalizing the main result in Z. Guo and J. Wei \cite{GW} to the case of $p \ne 2$ (with a simpler proof). In \cite{K} we called $w(t)$ {\em the generating solution}, and $w_0(t)$ {\em the guiding solution}.
\medskip

We apply a similar approach to a class of equations with polynomial $f(r,u)$ generalizing the well-known results of D.D. Joseph and T.S. Lundgren \cite{JL}, and to the $p$-Laplace version of the generalized Gelfand equation, where we easily recover the corresponding result of J. Jacobsen and K.  Schmitt \cite{JS}.
\smallskip

For each of the three classes of equations we show that along the solution curves (as $u(0) \ra \infty$), the solutions tend to a singular solution (for which $u(r) \ra \infty$, or $u'(r) \ra \infty$, as $r \ra 0$). Moreover, one can calculate the singular solutions explicitly, which is truly a remarkable feature of self-similar equations. Singular solutions were studied previously by many authors, including C. Budd and J.  Norbury \cite{BN}, F. Merle and L. A. Peletier \cite{MP}, and I. Flores \cite{F}.

\section{Parameterization of the solution curves}
\setcounter{equation}{0}
\setcounter{thm}{0}
\setcounter{lma}{0}

We begin with the $p$-Laplace version of the generalized Gelfand equation
\beq
\lbl{1}
\s\s\s \p (u')'+\frac{n-1}{r}\p (u') + \la \, r^{\alpha}e^u =0 \s \mb{for $0<r<1$,} \s u'(0)=0, \s u(1)=0 \,,
\eeq
where $\p(v)=v|v|^{p-2}$, $p>1$. Observe that $\p(sv)=s^{p-1} \p(v)$ for any constant $s>0$.
Assume that $u(0)=a>0$. We set $u=w+a$, $t=br$. The constants $a$ and $b$ are assumed to satisfy
\[
\la = b ^{\al +p}e^{-a} \,.
\]
Then (\ref{1}) becomes
\beq
\lbl{2}
\p (w')'+\frac{n-1}{t}\p (w')+  t^{\al}e^w=0, \s w(0)=0, \s w'(0)=0 \,.
\eeq
The solution of this problem $w(t)$, which  is a negative and decreasing function, is defined for all $t>0$, and it may be easily computed numerically. (Write this equation as $\left[ t^{n-1} \p (w')\right]'=-t^{n+\al -1} e^w<0$, and conclude that $t^{n-1} \p (w')<0$, and then $w'(t)<0$ for all $t$.) We have
\[
0=u(1)=a+w(b) \,,
\]
so that  $a=-w(b)$, and then  $\la = b ^{\al +p}e^{w(b)}$. The solution curve for (\ref{1}) is 
\[
(\la ,u(0))=\left(b ^{\al +p}e^{w(b)}\,, -w(b) \right) \,,
\]
parameterized by $b \in (0,\infty)$. The solution of (\ref{1}) at $b$ is $u(r)=w(br)-w(b)$.
It will be convenient 
to write the solution curve as 
\beq
\lbl{3}
(\la ,u(0))=\left(t ^{\al +p }e^{w(t)}\,, -w(t) \right) \,,
\eeq
parameterized by $t \in (0,\infty)$, and $w(t)$ is the solution of (\ref{2}). The solution of (\ref{1}) at the parameter value $t$ is $u(r)=w(tr)-w(t)$.
\medskip

We consider next the problem 
\beq
\lbl{4}
\p (u')'+\frac{n-1}{r}\p (u') + \la \, \frac{r^{\alpha}}{(1-u)^q} =0 \s\s \mb{for $0<r<1$} \,, 
\eeq
\[
u'(0)=u(1)=0 \,, \s 0<u(r)<1  \,,
\]
which arises in modeling of electrostatic micro-electromechanical systems (MEMS), see \cite{P}, \cite{GG}, \cite{GW}. Here $\la $ is a positive parameter, $q>0$ and $\al >0$ are constants, and as before $\p(v)=v|v|^{p-2}$, $p>1$. Any  solution $u(r)$ of (\ref{4}) is a positive and  decreasing function (by the maximum principle), so that $u(0)$ gives its maximum value. Our goal is to compute the solution curve $(\la,u(0))$. Let $1-u=v$. Then $v(r)$ satisfies
\beq
\lbl{5}
\p (v')'+\frac{n-1}{r}\p (v')= \la \, \frac{r^{\al}}{v^q}  \s\s \mb{for $0<r<1$,} \s v'(0)=0, \s v(1)=1 \,.
\eeq
Assume that $v(0)=a$. We scale $v(r)=aw(r)$, and $t=br$. The constants $a$ and $b$ are assumed to satisfy
\beq
\lbl{6}
\la =a^{p+q-1} b ^{\al +p} \,.
\eeq
Then (\ref{5}) becomes
\beq
\lbl{7}
\p (w')'+\frac{n-1}{t}\p (w')=  \frac{t^{\al}}{w^q}, \s w(0)=1, \s w'(0)=0 \,.
\eeq
The solution of this problem is a positive increasing function, which is defined for all $t>0$.
We have
\[
1=v(1)=aw(b) \,,
\]
and so $a=\frac{1}{w(b)}$, and then $\la =\frac{b ^{\al +p}}{w^{p+q-1}(b)}$. The solution curve $(\la,u(0))$ is
$
\left(\frac{b ^{\al +p}}{w^{p+q-1}(b)}\,, 1-\frac{1}{w(b)} \right)
$,
parameterized by $b \in (0,\infty)$.  It will be convenient 
to write the solution curve in the form
\beq
\lbl{8}
(\la,u(0))=\left(\frac{t ^{\al +p}}{w^{p+q-1}(t)}\,, 1-\frac{1}{w(t)} \right) \,,
\eeq
parameterized by $t \in (0,\infty)$.
In case $p=2$, this parameterization was first derived by J.A. Pelesko \cite{P},  and was then  used in  \cite{GG}. The solution of (\ref{4}) at $t$ is $u(r)=1-\frac{w(tr)}{w(t)}$.
\medskip

Finally, we consider the problem (with the constants $p>1$, $q>1$, $\al \geq 0$) 
\beq
\lbl{4a}
 \p (u')'+\frac{n-1}{r}\p (u') + \la \, r^{\alpha}(1+u)^q =0 \s \mb{for $0<r<1$}  \,,   
\eeq
\[
u'(0)=u(1)=0 \,,
\]
which was analyzed in case $p=2$ and $\alpha =0$ by D.D. Joseph and T.S. Lundgren \cite{JL}.
If we set  $1+u=v$, then $v(r)$ satisfies
\beq
\lbl{5a}
 \p (v')'+\frac{n-1}{r}\p (v')+ \la  r^{\al}v^q =0,  \s v'(0)=0, \s v(1)=1 \,.
\eeq
Assuming that $v(0)=a$, we scale $v(r)=aw(r)$, and $t=br$. The constants $a$ and $b$ are assumed to satisfy
\beq
\lbl{6a}
\la =\frac{ b ^{p+\al }}{a^{q-p+1}} \,.
\eeq
Then (\ref{5a}) becomes
\beq
\lbl{7a}
\p (w')'+\frac{n-1}{t}\p (w')+t^{\al}w^q=0, \s w(0)=1, \s w'(0)=0 \,.
\eeq
The solution of (\ref{7a}) satisfies $w'(t)<0$, so long as $w(t)>0$ (the function $t^{n-1} \p (w'(t))$ is zero at $t=0$, and its derivative is negative). It follows that either there is a $t_0$, so that $w(t_0)=0$ and $w(t)>0$ on $(0,t_0)$, or $w(t)>0$ on $(0,\infty)$ and $\lim _{t \ra \infty} w(t)=a \geq 0$. It is easy to see that $a=0$ in the second case. Indeed, assuming that $a>0$, we have $\left[t^{n-1} \p (w') \right]' \leq -a^q t^{n+\al -1}$, and integrating we conclude that $w(t) \leq 1-c t^{\gamma}$, with some $c>0$, and $\gamma=\frac{\al +p}{p-1}>0$, contradicting that $w(t)>0$ on $(0,\infty)$.

\begin{lma}\lbl{lma:6}
Assume that
\beq
\lbl{7b}
q>\frac{np-n+p+p\al}{n-p} \,.
\eeq
Then $w(t)>0$, and $w'(t)<0$ on $(0,\infty)$, with $\lim _{t \ra \infty} w(t)=0$.
\end{lma}

\pf
In view of the above remarks, we need to exclude the possibility that $w(t_0)=0$ and $w(t)>0$ on $(0,t_0)$.
Recall that for the equation
\[
\p (w')'+\frac{n-1}{t}\p (w')+f(t,w)=0 \,,
\]
the Pohozhaev function 
\[
P(t)=t^n \left[ (p-1)\p (w') w'+pF(t,w) \right]+(n-p)t^{n-1} \p (w')w
\]
is easily seen to satisfy
\[
P'(t)=t^{n-1} \left[ npF(t,w)-(n-p)wf(t,w)+pt F_t(t,w) \right] \,,
\]
where $F(t,w)=\int_0^w f(t,z) \, dz$, see e.g., \cite{K1}, p. $136$. Here 
\[
P'(t)=t^{n-1+\al}  \left[ \frac{np}{q+1}-(n-p)+ \frac{p \al }{q+1}\right] w^{q+1}<0 \,.
\]
Since $P(0)=0$, and $P(t_0)>0$, we have a contradiction.
\epf

As before, we have
\[
1=v(1)=aw(b) \,,
\]
and so $a=\frac{1}{w(b)}$, and then $\la =b ^{p+\al }w^{q-p+1}(b)$.  Under the condition (\ref{7b}), the solution curve $(\la,u(0))$ is
$
\left(b ^{p+\al }w^{q-p+1}(b)\,, \frac{1}{w(b)}-1 \right)
$,
parameterized by $b \in (0,\infty)$.  The solution at $b$ is $u(r)=\frac{w(br)}{w(b)}-1$. It will be convenient 
to write the solution curve in the form
\beq
\lbl{8b}
(\la,u(0))=\left(t ^{p+\al }w^{q-p+1}(t)\,, \frac{1}{w(t)}-1 \right) \,,
\eeq
parameterized by $t \in (0,\infty)$. The solution of (\ref{4a}) at $t$ is $u(r)=\frac{w(tr)}{w(t)}-1$.

\section{The equation modeling MEMS }
\setcounter{equation}{0}
\setcounter{thm}{0}
\setcounter{lma}{0}

We consider the problem (\ref{4}), whose solution curve is given by (\ref{8}), where $w(t)$ is the solution of (\ref{7}). We have $\la (t)=\frac{t ^{\al +p}}{w^{p+q-1}(t)}$, where $w(t)$ is the solution of (\ref{7}), and so
\[
\la '(t)=t^{\al +p-1}w^{-p-q} \left[(\al +p)w-t(p+q-1)w' \right] \,.
\]
We are interested in the roots of the function $(\al +p)w-t(p+q-1)w'$. If we set this function to zero
\[
(\al +p)w-t(p+q-1)w'=0 \,,
\]
then the general solution of this equation is 
\[
w(t)=c t^{\beta }, \s\s \beta=\frac{ \al +p}{p+q-1} \,.
\]
Quite remarkably, if we choose the constant $c=c_0=\left[\frac{1}{ \beta ^{p-1} \left[(p-1)(\beta-1)+n-1 \right]} \right]^{\frac{1}{p+q-1}}$, under the condition that
\beq
\lbl{*}
(p-1)(\beta-1)+n-1>0 \,,
\eeq
then 
\[
w_0(t)=c_0 t^{\beta }
\]
also solves the equation in (\ref{7}), along with $w(t)$. We shall  show that $w(t)$, the solution of the initial value problem (\ref{7}), tends to $w_0(t)$ as $t \ra \infty$, and the issue turns out to be whether $w(t)$ and $w_0(t)$ cross  infinitely many times as $t \ra \infty$.
\begin{lma}\lbl{lma:5}
Assume that $w(t)$ and $w_0(t)$ intersect infinitely many times. Then the solution curve of (\ref{4}) makes infinitely many turns.
\end{lma}

\pf
Assuming that $w(t)$ and $w_0(t)$ intersect infinitely many times, let $\{ t_n \}$ denote the points of intersection. At $\{ t_n \}$'s, $w(t)$ and $w_0(t)$ have different slopes (by uniqueness for initial value problems). Since $(\al +p)w_0(t_n)-t_n(p+q-1)w_0'(t_n)=0$, it follows that $(\al +p)w(t_n)-t_n(p+q-1)w'(t_n)<0$ ($>0$) if  $w(t)$ intersects $w_0(t)$ from below (above) at $t_n$. Hence, on any interval $(t_n,t_{n+1})$ there is a point $t_0$, where $(\al +p)w(t_0)-t_0(p+q-1)w'(t_0)=0$, i.e., $\la '(t_0)=0$, and $t_0$ gives a critical point. Since $\la '(t_n)$ and $\la '(t_{n+1})$ have different signs, the solution curve changes its direction over $(t_n,t_{n+1})$. 
\epf 

We shall need the following well-known Sturm-Picone's comparison theorem, see e.g., p. $5$ in  \cite{Kr}. 

\begin{lma}\lbl{lma:2} 
Let $u(t)$ and $v(t)$ be respectively classical solutions of 
\beq
\lbl{10a}
\left(a(t)u'\right)'+b(t)u=0 \,,
\eeq
\beq
\lbl{10b}
\left(a_1(t)v'\right)'+b_1(t)v=0 \,.
\eeq
Assume that the given differentiable functions $a(t)$, $a_1(t)$, and  continuous functions $b(t)$ and  $b_1(t)$, satisfy
\beq
\lbl{10c}
b_1(t) \geq b(t), \s \mbox{and} \s 0<a_1(t) \leq a(t) \s\s \mbox{for  $t \geq t_0>0$}.
\eeq
In case $a_1(t) = a(t)$ and $b_1(t) = b(t)$ for all $t$, assume additionally that $u(t)$ and $v(t)$ are not constant multiples of one another. Then, for  $t \geq t_0$,  $v(t)$ has a root between any two consecutive roots of $u(t)$.
\end{lma}

\begin{lma}\lbl{lma:3} 
Consider the equation
\beq
\lbl{12}
\s\s\s \left(a_0(t)\left(1+f(t) \right) v'\right)'+\frac{n-1}{t}a_0(t) \left(1+f(t) \right)v'+b_0(t)\left(1+g(t) \right)v=0 \,, 
\eeq
with given differentiable functions $a_0(t)>0$ and $f(t)$, and continuous functions   $b_0(t)>0$ and $g(t)$.
Assume that $\lim _{t \ra \infty} f(t)=\lim _{t \ra \infty} g(t)=0$, and there is an $\ep>0$ such that any solution of
\beq
\lbl{10d}
\left(a_0(t)\left(1+\ep \right) v'\right)'+\frac{n-1}{t}a_0(t) \left(1+\ep \right)v'+b_0(t)\left(1-\ep \right)v=0
\eeq
has infinitely many roots. Then any solution of (\ref{12}) has infinitely many roots.
\end{lma}

\pf
We rewrite (\ref{12}) in the form (\ref{10a}), with $a(t)=t^{n-1} a_0(t)\left(1+f(t) \right)$, and $b(t)=t^{n-1} b_0(t)\left(1+g(t) \right)$, and we rewrite (\ref{10d}) in the form (\ref{10b}), with $a_1(t)=t^{n-1} a_0(t)\left(1+\ep \right)$, and $b_1(t)=t^{n-1} b_0(t)\left(1-\ep \right)$.
For large $t$, the inequalities in (\ref{10c}) hold,  and the Lemma \ref{lma:2} applies.
\epf

The linearized equation for (\ref{7}) is
\[
\left( \p' (w')z' \right)'+\frac{n-1}{t}\p' (w')z'=-q t^{\al} w^{-q-1} z \,.
\]
At the solution $w=w_0(t)$, this becomes
\beq
\lbl{12a}
\left( a_0(t)z' \right)'+\frac{n-1}{t}a_0(t)z'+b_0(t) z =0\,,
\eeq
with $a_0(t)=\p' (w_0')=(p-1) c_0^{p-2} \beta ^{p-2} t^{(p-2)(\beta-1)}$, and $b_0(t)=q t^{\al} w_0^{-q-1}=qc_0^{-q-1} t^{\al- \beta (q+1)}$. One simplifies (\ref{12a}) to read 
\[
z''+\frac{\left[(p-2)(\beta-1)+n-1 \right]}{t}z'+\frac{q \beta  \left[(p-1)(\beta-1)+n-1 \right]}{(p-1) t^2}z=0 \,,
\]
which is  an Euler equation! The roots of its characteristic equation,
\[
r(r-1)+ \left[(p-2)(\beta-1)+n-1 \right]r+\frac{q \beta  \left[(p-1)(\beta-1)+n-1 \right]}{(p-1) } =0 \,,
\]
are complex valued, provided that 
\[
\left[(p-2)(\beta-1)+n-2 \right]^2< \frac{4q \beta  \left[(p-1)(\beta-1)+n-1 \right]}{p-1}\,.
\]
We write this inequality in the form
\beq
\lbl{20}
A\beta ^2+B\beta-C>0 \,,
\eeq
with $A=4(p-1)q-(p-1)(p-2)^2$, $B=4q(n-p)-2 (p-1)(p-2)(n-p)$, and $C=(p-1)(n-p)^2$. We shall have $A>0$, provided that
\beq
\lbl{20.1}
4q-(p-2)^2>0 \,.
\eeq
For (\ref{20}) to hold, we need $\beta=\frac{ \al +p}{p+q-1}$ to be greater than the larger root of this quadratic, i.e., $\beta >\frac{-B+\sqrt{B^2+4AC}}{2A}$ (assuming (\ref{20.1})), which gives
\beq
\lbl{21}
\s \; \frac{ \al +p}{p+q-1}>\frac{(p-n) \left(2q-p^2+3p-2 \right)+2 |n-p| \sqrt{q(p+q-1)}}{(p-1) \left[4q-(p-2)^2 \right]}  \,.
\eeq
\begin{thm}\lbl{thm:2}
Assume that $q>0$, $p >1$, with
\beq
\lbl{21.1}
(p-1)(\beta-1)+n-1>\beta \,,
\eeq
and the conditions  (\ref{20.1}) and (\ref{21}) hold. Then the solution curve of (\ref{4}) makes infinitely many turns. Moreover, along this curve (as $u(0) \ra \infty$), $\la \ra \la _0=\frac{1}{c_0^{q-1}}= \beta ^{p-1} \left[(p-1)(\beta-1)+n-1 \right]$, and $u(r) $ tends to $1- r^{\beta}$ for $r \ne 0$, which is a  solution of the equation in (\ref{4}).  
\end{thm}

\pf
In view of Lemma \ref{lma:5}, we need to show that $w(t)$ and $w_0(t)$ intersect infinitely many times. Let $P(t)=w(t)-w_0(t)$. Then $P(t)$ satisfies
\beq
\lbl{23}
\left(a(t)P' \right)'+\frac{n-1}{t}a(t) P' +b(t)P=0 \,,
\eeq
where
\beq
\lbl{24}
a(t)=\int_0^1 \p ' \left(sw'(t)+(1-s)w_0'(t) \right) \, ds \,,
\eeq
\beq
\lbl{25}
b(t)=q \, t^{\al } \int_0^1 \frac{1}{\left[s w(t)+(1-s) w_0(t) \right]^{q+1}} \, ds \,.
\eeq
We claim that it is impossible for $P(t)$ to keep the same sign over some infinite interval $(t_0, \infty)$, and tend to a constant as $t \ra \infty$.
Assuming the contrary, write
\[
a(t)=(p-1) \left( w_0' \right)^{p-2}\int_0^1   \left|s \frac{w'(t)}{w'_0(t)} +(1-s)   \right|^{p-2} \, ds=a_0(t)\left(1+o(1) \right) \,,
\]
\[
b(t)=q \, t^{\al }\frac{1}{w_0^{q+1}}\int_0^1  \frac{1}{ \left[s \frac{w(t)}{w_0(t)} +(1-s)   \right]^{q+1}} \, ds=b_0(t)\left(1+o(1) \right) \,.
\]
as $t \ra \infty$. (Observe that $\frac{w(t)}{w_0(t)} \ra 1$, since $P(t)$ tends to a constant, and $\frac{w'(t)}{w'_0(t)} \ra 1$, by L'Hospital's rule, as $t \ra \infty$.)
Since Euler's equation (\ref{12a}) has infinitely many roots on $(t_0,\infty)$, we conclude by Lemma \ref{lma:3} that $P(t)$ must vanish on that interval too, a contradiction.
\medskip

Next we show that if $P(t_0)=0$, then $P(t)$ remains bounded for all $t>t_0$. Assume that $P'(t_0)<0$, and the case when $P'(t_0)>0$ is similar. Then  $P(t)<0$ for $t>t_0$, with $t-t_0$ small. From (\ref{23}), $t^{n-1}a(t)P'(t)$ is increasing  for $t>t_0$,  so that
\[
P'(t)>-\frac{a_0}{a(t)t^{n-1}}, \s \mbox{for $t>t_0$}\s\s  (\mbox{with}\; a_0=-t_0^{n-1}a(t_0)P'(t_0)>0) \,.
\]
Since solutions of the linear equation (\ref{23}) cannot go to infinity over a bounded interval, we may assume that $t_0$ is large, and then by the above
$a(t) \sim a_0(t) \sim a_1 t^{(p-2)(\beta -1)}$ for $t>t_0$, and some $a_1 >0$. It follows that for some $a_2 >0$
\beq
\lbl{25.1}
P'(t)>-\frac{a_2}{t^{n-1+(p-2)(\beta -1)}}=-\frac{a_2}{t^{1+\epsilon}}, \s \mbox{for $t>t_0$} \,,
\eeq
with $\epsilon=n-2+(p-2)(\beta-1)>0$, in view of (\ref{21.1}).
Integrating over $(t_0,t)$, and using that $n \geq 3$, we conclude the boundness of $P(t)$, so long as $P(t)<0$. If another root of $P(t)$ is encountered, we repeat the argument. Hence, $P(t)$ remains bounded for all $t>t_0$.
\medskip

From the equation (\ref{23}), we see that $P(t)$ cannot have points of positive minimum or points of negative maximum. We claim that if $P(t)$ has one root, it has infinitely many roots. Indeed, assume that $P(t_1)=0$, and say  $P'(t_1)>0$. For $t>t_1$, $P(t)$ remains bounded, but cannot tend to a constant. Hence, $P(t)$ will have to turn back and become decreasing, but it cannot have a positive local minimum, or tend to a constant. Hence, $P(t_2)=0$ at some 
$t_2>t_1$, and so on. 
\medskip

We have $P(0)=1$, so that $\left(t^{n-1}a(t)P'(t) \right)'<0$ for small $t>0$. The function $q(t) \equiv t^{n-1}a(t)P'(t)$ satisfies $q(0)=0$ and $q'(t)<0$, and so $q(t)<0$. It follows that $P'(t)<0$ for small $t>0$. Since $P(t)$ cannot turn around, or tend to a constant, we conclude the existence of the first root $t_1$ of $P(t)$, implying the existence of infinitely many roots.
\medskip

We show next that $w(t) \ra w_0(t)$ as $t \ra \infty$. 
Let $t_k$ and $t_{k+1}$ be two consecutive roots of $P(t)$, and $P'(t_k)<0$, so that $P(t)<0$ on $(t_k,t_{k+1})$. Let $\tau _k$ be the unique point of minimum of $P(t)$ on $(t_k,t_{k+1})$.
For negative $P(t)$ we have the inequality (\ref{25.1}), with $t_k$ in place of $t_0$.
Integrating this  inequality over $(t_k,\tau _k)$,  we get 
\[
P(\tau _k) >\bar c \left(\tau _k ^{-\epsilon}-t _k ^{-\epsilon} \right)  \s (\mbox{ with some $\bar c>0$}) \,,
\]
which implies that $|P(\tau _k)| \ra 0$, as $k \ra \infty$. The case when $P'(t_k)>0$ is similar, so that $w(t) \ra w_0(t)$ along the solution curve.
Since $u(r)=1-\frac{w(tr)}{w(t)}$, it follows that along the solution curve
$u(r)$ tends to $1-\frac{w_0(tr)}{w_0(t)}=1-r^{\beta}$, while $\la(t)$ tends to $\frac{1}{c_0^{q-1}}$.
\epf

Observe that in case $\beta \in (0,1)$, the limiting solution $1-r^{\beta}$ is {\em singular}, because $u'(0)$ is not defined. Notice also that the condition (\ref{21.1}) implies (\ref{lma:5}). Finally, observe that in case $\beta \in (0,1)$ the condition (\ref{21.1}) implies that $n \geq 2$. Indeed, we can rewrite (\ref{21.1})
as $n > 2 \beta+p(1-\beta)$, which is a point between $p>1$, and $2$.
\medskip

One special case when this theorem applies is the following. Assume that $n \geq p$, so that (\ref{21}) becomes
\[
\frac{ \al +p}{p+q-1}>(n-p) \, \frac{2  \sqrt{q(p+q-1)} +p^2-3p+2 -2q}{(p-1) \left[4q-(p-2)^2 \right]}  \,.
\]
Then (\ref{21}) holds, provided that
\beq
\lbl{25.2}
2  \sqrt{q(p+q-1)} +p^2-3p+2 -2q > 0 \,,
\eeq
\[
4q>(p-2)^2 \,,
\]
\[
p \leq n < p+\frac{ (\al +p)(p-1) \left[4q-(p-2)^2 \right]}{(p+q-1) \left(2  \sqrt{q(p+q-1)} +p^2-3p+2 -2q \right) } \,.
\]
Observe that the third inequality ($n \geq p$) implies that the condition (\ref{*}) holds, and the second inequality is just (\ref{20.1}).
Hence, the three inequalties in (\ref{25.2}) imply the theorem. 
In case $p=2$, the first and the second inequalities hold automatically, while  the third one  gives the condition in Z. Guo and J. Wei \cite{GW}.

\section{The generalized Joseph-Lundgren problem}
\setcounter{equation}{0}
\setcounter{thm}{0}
\setcounter{lma}{0}

We now study the problem (\ref{4a}). Its solution curve is represented by (\ref{8b}), under the condition (\ref{7b}), where $w(t)$ is the solution of (\ref{7a}).
In particular, $\la (t)=t ^{p+\al }w^{q-p+1}(t)$, and we wish to know  how many times this function changes the direction of monotonicity for $t \in (0,\infty)$. (Here $w(t)$ is the  solution of (\ref{7a}), the generating solution.)
Compute
\[
\la '(t)=t^{p+\al -1} w^{q-p}(t) \left[(p+\al)w(t) +(q-p+1)tw'(t) \right] \,,
\]
so that we are interested in the roots of the function $(p+\al)w +(q-p+1)tw'$. If we set this function to zero
\[
(p+\al)w +(q-p+1)tw'=0 \,,
\]
then the general  solution of this equation is $w(t)=at^{-\beta}$, with $\beta=\frac{p+\al}{q-p+1}$.
If we choose the constant $a$ as
\[
a=a_0=\left[(n-p)\beta^{p-1}-(p-1)\beta^p \right]^{\frac{1}{q-p+1}}
\]
then $w_0(t)=a_0t^{-\beta}$ is a solution of (\ref{7a}), the guiding solution (we have $(n-p)\beta^{p-1}-(p-1)\beta^p>0$, under the condition (\ref{7b}), if $n>p$).
\begin{lma}\lbl{lma:15}
Assume that $w(t)$ and $w_0(t)$ intersect infinitely many times. Then the solution curve of (\ref{4a}) makes infinitely many turns.
\end{lma}

\pf
Indeed, assuming that $w(t)$ and $w_0(t)$ intersect infinitely many times, let $\{ t_n \}$ denote their points of intersection. At $\{ t_n \}$'s, $w(t)$ and $w_0(t)$ have different slopes (by uniqueness for initial value problems). Since $(p+\al)w_0(t_n) +(q-p+1)t_nw_0'(t_n)=0$, it follows that $(p+\al)w(t_n) +(q-p+1)t_nw'(t_n)>0$ ($<0$) if  $w(t)$ intersects $w_0(t)$ from below (above) at $t_n$. Hence, on any interval $(t_n,t_{n+1})$ there is a point $t_0$, where $(p+\al)w(t_0) +(q-p+1)t_0w'(t_0)=0$, i.e., $\la '(t_0)=0$, and $t_0$ is a critical point. Since $\la '(t_n)$ and $\la '(t_{n+1})$ have different signs, the solution curve changes its direction over $(t_n,t_{n+1})$. 
\epf 

\medskip

The linearized equation for (\ref{7a}) is
\[
\left( \p' (w')z' \right)'+\frac{n-1}{t}\p' (w')z'+q t^{\al} w^{q-1} z=0 \,.
\]
At the solution $w=w_0(t)$, this becomes
\beq
\lbl{30}
\left( a_0(t)z' \right)'+\frac{n-1}{t}a_0(t)z'+b_0(t) z =0\,,
\eeq
with $a_0(t)=\p' (w_0')$, and $b_0(t)=q t^{\al} w_0^{q-1}$.
One simplifies (\ref{30}) to Euler's equation 
\beq
\lbl{31}
z''+\frac{\left[-(\beta+1)(p-2)+n-1 \right]}{t}z'+\frac{q a_0^{q-p+1}}{(p-1) {\beta}^{p-2} t^2}z=0 \,.
\eeq

Let us consider  first the case when $p=2$ and $\al =0$, and $n>2$. Then $\beta =\frac{2}{q-1}$, $a_0=\left[\beta(n-\beta-2) \right]^{\frac{1}{q-1}}$, and the equation (\ref{31}) becomes
\[
t^2z''+(n-1)tz'+q \beta (n-\beta-2)z=0\,.
\]
Its characteristic equation 
\[
r(r-1)+(n-1)r+q \beta (n-\beta-2)=0
\]
has the roots
\[
r=\frac{-(n-2) \pm \sqrt{(n-2)^2-4q\beta (n-\beta-2)}}{2} \,.
\]
These roots are complex if 
\[
(n-2)^2-4q\beta (n-2)+4q\beta^2 <0 \,.
\]
On the left we have a quadratic in $n-2$, with two positive roots. The largest value of $n-2$, for which this inequality holds, corresponds to the larger root of this quadratic, i.e.,
\beq
\lbl{32}
n-2 < \frac{4q}{q-1}+4 \sqrt{ \frac{q}{q-1}} \,.
\eeq
We shall show that infinitely many solutions occur if (\ref{32}) holds, and 
\beq
\lbl{32.1}
q>\frac{n+2}{n-2} \,.
\eeq
(The last condition ensures  that the generating solution $w(t)$ is defined for all $t>0$, by Lemma \ref{lma:6}.)  In terms of $n$, the conditions (\ref{32}) and (\ref{32.1}) imply
\beq
\lbl{32.2}
\frac{2+2q}{q-1}<n<2+\frac{4q}{q-1}+4 \sqrt{ \frac{q}{q-1}} \,,
\eeq
which is the condition from \cite{JL} (it implies that $n>2$). Thus we shall recover the following classical theorem of D.D. Joseph and T.S. Lundgren \cite{JL}.

\begin{thm}\lbl{thm:20}
Assume that  the conditions  (\ref{32}) and  (\ref{32.1}) hold (or (\ref{32.2}) holds). 
Then the solution curve of (\ref{4a}) makes infinitely many turns. Moreover, along this curve (as $u(0) \ra \infty$), $\la \ra \la _0=a_0^{q-1}$, and $u(r) $ tends to $ r^{-\beta}-1$ for $r \ne 0$, which is a singular solution of the equation in (\ref{4a}).  
\end{thm}
We shall give a proof of more general result below.
\medskip

For general $p$ and $\al$, the characteristic equation for (\ref{31}) is
\beq
\lbl{32.3}
r(r-1)+A r+B=0 \,,
\eeq
with $A=-\beta (p-2)+n-p+1$, and $B=\frac{q(n-p)}{p-1} \beta-q \beta^2$. The roots of (\ref{32.3})
\[
r=\frac{-(A-1) \pm \sqrt{(A-1)^2-4B}}{2}
\]
are complex, provided that
\[
(A-1)^2-4B<0 \,,
\]
which simplifies to
\beq
\lbl{32.4}
(n-p)^2-\theta (n-p)+\gamma<0 \,,
\eeq
with
\beq
\lbl{32.5}
\theta = 2\beta(p-2)+\frac{4q\beta}{p-1} \,, \s \gamma = (p-2)^2\beta^2+4q\beta^2\,.
\eeq
On the left in (\ref{32.4}) we have a quadratic in $n-p$, with two positive roots. The largest value of $n-p$, for which the inequality (\ref{32.4}) holds, corresponds to the larger root of this quadratic, i.e.,
\beq
\lbl{32.6}
n-p < \frac{\theta+\sqrt{\theta^2-4\gamma}}{2} \,.
\eeq
We shall show that infinitely many solutions occur if the conditions (\ref{7b}) and (\ref{32.6}) hold. In terms of $n$, the conditions (\ref{7b}) and (\ref{32.6}) imply that
\beq
\lbl{32.7}
\frac{pq+p+p\al}{q-p+1}<n < p+\frac{\theta+\sqrt{\theta^2-4\gamma}}{2} \,.
\eeq
The first inequality in (\ref{32.7}) implies that 
\beq
\lbl{32.03}
(\beta+1)(p-2)<n-2 \,,
\eeq
which in turn gives that $n>p$.
\begin{thm}\lbl{thm:21}
Assume that  the conditions  (\ref{7b}) and  (\ref{32.6}) hold (or (\ref{32.7}) holds). 
Then the solution curve of (\ref{4a}) makes infinitely many turns. Moreover, along this curve (as $u(0) \ra \infty$), $\la \ra \la _0=a_0^{q-1}$, and $u(r) $ tends to $ r^{-\beta}-1$ for $r \ne 0$, which is a singular solution of the equation in (\ref{4a}).  
\end{thm}

\pf
In view of Lemma \ref{lma:15}, we need to show that $w(t)$ and $w_0(t)$ intersect infinitely many times, and they tend to each other as $t \ra \infty$. Let $P(t)=w(t)-w_0(t)$. Then $P(t)$ satisfies
\beq
\lbl{53}
\left(a(t)P' \right)'+\frac{n-1}{t}a(t) P' +b(t)P=0 \,,
\eeq
where
\beq
\lbl{54}
a(t)=\int_0^1 \p ' \left(sw'(t)+(1-s)w_0'(t) \right) \, ds \,,
\eeq
\beq
\lbl{55}
b(t)=q t^{\al} \,  \int_0^1 \left[s w(t)+(1-s) w_0(t) \right]^{q-1} \, ds \,.
\eeq
Since both $w(t)$ and $w_0(t)$ tend to zero as $t \ra \infty$ (see Lemma \ref{lma:6}), we conclude that $P(t) \ra 0$ as $t \ra \infty$. This simplifies the proof considerably, compared to Theorem \ref{thm:2}.
We claim that it is impossible for $P(t)$ to keep the same sign over some infinite interval $(t_0, \infty)$.
Assuming the contrary, write ($a_0(t)$ and $b_0(t)$ were defined in (\ref{30}))
\[
a(t)=(p-1) \left( -w_0' \right)^{p-2}\int_0^1   \left|s \frac{w'(t)}{w'_0(t)} +(1-s)   \right|^{p-2} \, ds=a_0(t)\left(1+o(1) \right) \,,
\]
\[
b(t)=q t^{\al} \, w_0^{q-1}\int_0^1  \left[s \frac{w(t)}{w_0(t)} +(1-s)   \right]^{q-1} \, ds=b_0(t)\left(1+o(1) \right) \,.
\]
as $t \ra \infty$. (Observe that $\frac{w(t)}{w_0(t)} \ra 1$, since $P(t)$ tends to a constant, and $\frac{w'(t)}{w'_0(t)} \ra 1$, by L'Hospital's rule, as $t \ra \infty$.)
Since Euler's equation (\ref{12a}) has infinitely many solutions on $(t_k,\infty)$, we conclude by Lemma \ref{lma:3} that $P(t)$ must vanish on that interval too, a contradiction. It follows that $P(t)$ has infinitely many roots, which implies that $w(t)$ and $w_0(t)$ have infinitely many points of intersection, and hence the solution curve makes infinitely many turns.
\medskip

Since $u(r)=\frac{w(tr)}{w(t)}-1$, it follows that along the solution curve
$u(r)$ tends to $\frac{w_0(tr)}{w_0(t)}-1=r^{-\beta}-1$ for $r \ne 0$.
\epf

\section{The generalized Gelfand problem}
\setcounter{equation}{0}
\setcounter{thm}{0}
\setcounter{lma}{0}

We now use the representation (\ref{3}) for the solution curve of (\ref{1}). In particular, $\la (t)=t ^{\al +p}e^{w(t)}$, where $w(t)$ is the solution of (\ref{2}), and the issue is how many times this function changes its direction of monotonicity for $t \in (0,\infty)$. 
Compute
\[
\la '(t)=t e^{w} \left(\al +p+tw' \right) \,,
\]
so that we are interested in the roots of the function $\al +p+tw'$. If we set this function to zero
\[
\al +p+tw'=0 \,,
\]
then the solution of this equation is of course $w(t)=a-(\al +p) \ln t$. Quite surprisingly, if we choose the constant $a=a_0=\ln \left[ (n-p)(\al +p)^{p-1} \right]$, assuming that $n >p$, then 
\[
w_0(t)=\ln \left[ (n-p)(\al +p)^{p-1} \right]-(\al +p) \ln t
\]
is a solution of the equation in (\ref{2})!
We shall  show that $w(t)$ (the solution of the initial value problem (\ref{2})) tends to $w_0(t)$ as $t \ra \infty$, and give a condition for  $w(t)$ and $w_0(t)$ to cross  infinitely many times as $t \ra \infty$.
\begin{lma}\lbl{lma:1}
Assume that $w(t)$ and $w_0(t)$ intersect infinitely many times. Then the solution curve of (\ref{1}) makes infinitely many turns.
\end{lma}

\pf
Indeed, assuming that $w(t)$ and $w_0(t)$ intersect infinitely many times, let $\{ t_n \}$ denote the points of intersection. At $\{ t_n \}$'s, $w(t)$ and $w_0(t)$ have different slopes (by uniqueness for initial value problems). Since $\al +p+t_nw'_0(t_n)=0$, it follows that $\al +p+t_nw'(t_n)>0$ ($<0$) if  $w(t)$ intersects $w_0(t)$ from below (above) at $t_n$. Hence, on any interval $(t_n,t_{n+1})$ there is a point $t_0$, where $\al +p+t_0w'(t _0)=0$, i.e., $\la '(t_0)=0$, and $t_0$ is a critical point. Since $\la '(t_n)$ and $\la '(t_{n+1})$ have different signs, the solution curve changes its direction over $(t_n,t_{n+1})$. 
\epf 

The linearized equation for (\ref{2}) is
\[
\left( \p' (w')z' \right)'+\frac{n-1}{t}\p' (w')z'+t^{\al}e^wz=0 \,.
\]
At the solution $w=w_0(t)$, this becomes
\beq
\lbl{9}
 \left( a_0(t)z' \right)'+\frac{n-1}{t}a_0(t)z'+b_0(t)z=0\,,
\eeq
with $a_0(t)=\p' (w_0')=\frac{(p-1)(p+\al)^{p-2}}{t^{p-2}}$, and $b_0(t)=t^{\al}e^{w_0}=\frac{(n-p)(p+\al)^{p-1}}{t^{p}}$.
Simplifying (\ref{9}) gives
\[
(p-1) t^2z''+(p-1)(n-p+1)tz'+(n-p)( p+\al)z=0 \,,
\]
which is   Euler's equation! Its characteristic equation 
\[
(p-1)\, r(r-1)+(p-1)(n-p+1)\, r+(n-p)( p+\al)=0
\]
has the roots
\[
r=\frac{-(p-1)(n-p) \pm \sqrt{((p-1)(n-p) \left[p-1)(n-p)-4(p+\al) \right]}}{2(p-1)} \,.
\]
The roots are complex if $n-p>0$, and the quantity in the square brackets is negative (the opposite inequalities lead to a vacuous condition), i.e., when
\beq
\lbl{10}
p<n<\frac{p^2+3p+4 \al}{p-1} \,.
\eeq

We now easily recover the following result of J. Jacobsen and K.  Schmitt \cite{JS}, which was a generalization of the famous theorem of D.D. Joseph and T.S. Lundgren \cite{JL}.
\begin{thm}\lbl{thm:1}
Assume that the condition (\ref{10}) holds. Then the solution curve of (\ref{1}) makes infinitely many turns. Moreover, along this curve (as $u(0) \ra \infty$), $\la \ra e^{a _0}=(n-p)(p+\al)^{p-1}$, and $u(r) $ tends to $-(p+\al) \ln r$ for $r \ne 0$, which is a  singular solution of the equation in (\ref{1}).
\end{thm}

\pf
We follow the proof of the Theorem \ref{thm:2}.
In view of Lemma \ref{lma:1}, we need to show that $w(t)$ and $w_0(t)$ intersect infinitely many times. Let $P(t)=w(t)-w_0(t)$. Then $P(t)$ satisfies
\beq
\lbl{14}
\left(a(t)P' \right)'+\frac{n-1}{t}a(t) P' +b(t)P=0 \,,
\eeq
where
\beq
\lbl{16}
a(t)=\int_0^1 \p ' \left(sw'(t)+(1-s)w_0'(t) \right) \, ds \,,
\eeq
\beq
\lbl{17}
b(t)=t^{\al} \int_0^1 e^{sw(t)+(1-s)w_0(t)} \, ds \,.
\eeq
Compared with the proof of the Theorem  \ref{thm:2}, we have a complication here: in case $P(t)$ tends to a constant $p_0$ as $t \ra \infty$, we cannot conclude that $b(t)=b_0 (t)(1+o(1))$, unless $p_0=0$.
\medskip

We claim that it is impossible for $P(t)$ to keep the same sign over some infinite interval $(t_0, \infty)$, and tend to a constant $p_0 \ne 0$  as $t \ra \infty$. Assume, on the contrary, that $P(t)>0$ on $(t_0, \infty)$, and $\lim _{t \ra \infty} P(t)=p_0 > 0$. We may assume that 
\beq
\lbl{17b}
P(t)> \frac12 p_0>0  \s \mbox{ on $(t_1, \infty)$, with some $t_1>t_0$} \,.
\eeq
Write (\ref{14}) as 
\beq
\lbl{19a}
\left( t^{n-1} a(t)P' \right)'=-t^{n-1} b(t)P \,.
\eeq
As before,
\beq
\lbl{20a}
a(t)=a_0(t) \left(1+f(t) \right), \s\s \mbox{with $f(t) \ra 0\;$ as $t \ra \infty$} \,.
\eeq
Writing $b(t)=t^{\al}e^{w_0(t)} \int_0^1 e^{s P(t)} \, ds$, we see that 
\beq
\lbl{21a}
b(t)=b_0(t) \left(p_1+g(t) \right) \,,
\eeq
with $p_1=\int_0^1 e^{s p_0} \, ds>1$, and $g(t) \ra 0$ as $t \ra \infty$. By (\ref{19a}), (\ref{17b}),  and (\ref{21a})
\[
\left( t^{n-1} a(t)P' \right)'<-c_1t^{n-p-1} \s \mbox{on $(t_1,\infty)$}\,,
\]
for some constant $c_1>0$. Integrating this inequality over $(t_1,t)$, we get
\beq
\lbl{22a}
t^{n-1} a(t)P'<c_2-c_3t^{n-p} \s \mbox{on $(t_1,\infty)$}\,,
\eeq
for some constants $c_2>0$,  and $c_3>0$ (using that $n>p$). By (\ref{20a})
\[
a(t) >c_4t^{-p+2} \s \mbox{on $(t_2,\infty)$}\,,
\]
for some constants $c_4>0$,  and $t_2>t_1$. Using this in (\ref{22a}), we have
\[
P'<\frac{c_2}{c_4} t^{-n+p-1}-\frac{c_3}{c_4} \, t^{-1} \s \mbox{on $(t_2,\infty)$}\,.
\]
Integrating this over $(t_2,t)$, and using that $n>p$
\[
P(t)<c_5+\frac{c_2}{c_4(-n+p)} t^{-n+p}-\frac{c_3}{c_4} \ln t <c_5-\frac{c_3}{c_4} \ln t \,,
\]
for some constant $c_5>0$. Hence, $P(t)$ has to vanish at some $t>t_2$, contradicting the assumption that $P(t)>0$ on $(t_0,\infty)$. This proves that 
$p_0=0$.
We conclude that $p_1=1$ in (\ref{21a}), and the rest of the proof is similar to that of Theorem  \ref{thm:2}. 
\epf

If $p=2$ and $\al=0$, the condition (\ref{10}) becomes $2<n<10$, the classical condition of D.D. Joseph and T.S. Lundgren \cite{JL}.

\end{document}